\documentclass{amsart}

\usepackage{amssymb,latexsym}

\theoremstyle{plain}

\theoremstyle{definition}

\theoremstyle{remark}

\numberwithin{equation}{section}

\newcommand{\otter}{\textsc{Otter}}
\newcommand{\sem}{\textsc{Sem}}
\newcommand{\mace}{\textsc{Mace}}

\newcommand{\M}{\negthinspace\wedge\negthinspace}
\newcommand{\J}{\negthinspace\vee\negthinspace}

\begin{document} %

\title{Yet Another Single Law for Lattices}

\author{William McCune}
\address{Mathematics and Computer Science Division\\
        Argonne National Laboratory\\
        Argonne, Illinois  60419\\
        U.S.A.}
\email{mccune@mcs.anl.gov}
\urladdr{http://www.mcs.anl.gov/\~{}mccune}
\thanks{The research of the first author was
        supported by the Mathematical, Information, and
        Computational Sciences Division subprogram of the Office of Advanced
        Scientific Computing Research, U.S. Department of Energy, under
        Contract W-31-109-Eng-38.}

\author{R. Padmanabhan}
\address{
   Department of Mathematics\\
   University of Manitoba\\
   Winnipeg, Manitoba  R3T 2N2\\
   Canada}
\email{padman@cc.umanitoba.ca}
\urladdr{http://home.cc.umanitoba.ca/\~{}padman}
\thanks{The research of the second author was supported by an
         operating grant from NSERC of Canada (OGP8215).}

\author{Robert Veroff}
\address{
   Department of Computer Science\\
   University of New Mexico\\
   Albuquerque, New Mexico  87131\\
   U.S.A.}
\email{veroff@cs.unm.edu}
\urladdr{http://www.cs.unm.edu/\~{}veroff}



\date{}

\begin{abstract}
In this note we show that the equational theory of all 
lattices is defined by the single absorption law
\[
(((y\J x)\M x)\J (((z\M (x\J x))\J (u\M x))\M v))\M (w\J ((s\J x)\M (x\J t))) = x.
\]
This identity of length 29 with 8 variables is shorter than 
previously known such equations defining lattices.
\end{abstract}

\maketitle

\section{Introduction} \label{sec-intro}

Given a finitely based equational theory of algebras, it is natural to
determine the least number of equations needed to define that theory.
Researchers have known for a long time that all finitely based
group theories are one based \cite{higman-neumann}.
Because of lack of a cancellation law in lattices (i.e., the
absence of some kind of a subtraction operation), it was widely
believed that the equational theory of lattices cannot be defined
by a single identity.  This belief was further strengthened by the fact that
two closely related varieties, semilattices and distributive lattices,
were shown to be not one based \cite{potts:1965,mckenzie-lattice}.

In the late 1960s researchers attempted to formally prove that lattice
theory is not one based by trying to show that no set of absorption
laws valid in lattices can capture associativity.
Given one such attempt \cite{mckenzie-to-rp:1968},
Padmanabhan pointed out that Sholander's 2-basis for
distributive lattices \cite{sholander:lt} caused the
method to collapse.
This failure of the method
led to a proof of existence of a single
identity for lattices \cite{mckenzie-lattice}.  For a
partial history of various single identities defining lattices, their
respective lengths, and so forth, see the latest book on lattice theory by
G.~Gr\"atzer \cite[p. 477]{gratzer:1998}.

In this paper we present an absorption law of length 29 with 8
variables that characterizes the equational theory of all lattices.
To the best of our knowledge, the shortest previously known
single identity has length 77 with 8 variables \cite{veroff:web-lt77}.
Table~\ref{summary-tab} summarizes the results.
\begin{table}[ht] \label{summary-tab} \small
\caption{Discovery of Single Identities for Lattice Theory}
\begin{tabular}{llrr}
Authors                    & Reference  & Length  & Variables \\
\hline                    
R.~McKenzie                & 1970 \cite{mckenzie-lattice} & 300,000 &  34 \\
R. Padmanabhan             & 1977 \cite{rp-major}         & 243     &   7 \\
W. McCune, R. Padmanabhan  & 1996 \cite{wwm-rp:lt-wal}    & 79      &   7 \\
R. Veroff                  & 2001 \cite{veroff:web-lt77}  & 77      &   8 \\
W. McCune, R. Padmanabhan, R. Veroff 
                           & 2002                         & 29      &   8 \\
\hline
\end{tabular}
\end{table}

\section{Methodology}

The previously known single identities for lattice theory
were found by procedures that take a basis in the
form of absorption equations and reduce the size of the basis
to a single equation.  Such procedures typically produce
very large equations.  The single identities presented here
were found by enumerating lattice identities,
filtering them through sets of nonlattices, then
trying to prove automatically that the surviving equations are
single identities.  Several programs were used.
\begin{description}
\item[Eq-enum] enumerates equations of the form $\alpha = x$,
where $\alpha$ is in terms of meet, join, and variables.
Each variable in $\alpha$, aside from $x$, has one occurrence
(most-general absorption lattice identities have this property),
and neither the leftmost nor the rightmost variable in $\alpha$
is $x$ (such identities are eliminated by projection models).
\item[Lattice-filter] takes a stream of equations, uses Whitman's algorithm
to decide which are lattice identities, and discards the nonidentities.
\item[S{\small EM} \cite{sem} and M{\small ACE} \cite{mace2}]
search for small, finite nonlattice models of equations.
\item[Model-filter] takes a set of finite structures (nonlattices in
our case) and a stream of equations and discards equations that are
true in any of the structures.  Most of the nonlattices were found
by \sem\ or \mace; several were constructed by hand while examining
candidate identities.
\item[O{\small TTER} \cite{otter3}] searches for proofs of first order
and equational statements.  In this case, it took candidate identities
and tried for several seconds to prove basic lattice properties
such as commutativity, associativity, idempotence, and absorption identities.
\end{description}
The general method was to apply the preceding programs
and incrementally build a set of nonlattice structures
by using \sem\ and \mace\ to search for nonlattice models
of the current candidate.  Some of these structures were then
added to the set and used for filtering the subsequent candidates.

The programs were combined into a single program that was driven
by the choice of which equations to enumerate (for example length
25 with 7 variables), and which nonlattices to use for filtering.
It was run on several hundred processors, usually in jobs of 10--20
hours, over a period of several weeks.
If the computation had been done on one processor, it
would have taken several years.  About half a trillion equations
were enumerated, and the set of nonlattices grew to be several
thousand, most of size 4.  About 100,000 candidates survived the
model filter program, and \otter\ proved basic lattice properties
for several hundred of those.  Further \otter\ searches
on those candidates showed the following two to be single identities.
\begin{equation}
(((y\J x)\M x)\J (((z\M (x\J x))\J (u\M x))\M v))\M (w\J ((s\J x)\M (x\J t))) = x
\tag{A1}
\end{equation}
\begin{equation}
(((y\J x)\M x)\J (((z\M (x\J x))\J (u\M x))\M v))\M (((w\J x)\M (s\J x))\J t) = x
\tag{A2}
\end{equation}
The first proof that (A1) is a basis took \otter\ several
minutes and was more than 250 steps.
The standard lattice theory
6-basis (commutativity, associativity, absorption) was derived.
We then had \otter\ prove McKenzie's 4-basis (given in the next section),
which produced a proof of about 170 steps in less than one minute.
The 50-step proof given below was produced by Larry Wos, who used various
methods to simplify \otter 's proof of the McKenzie basis.

The search for single identities was not complete;
that is, many shorter equations and equations with fewer variables
were considered for which we could find neither countermodels nor proofs.
Therefore, whether there exists a shorter single identity is an open
question.

\section{Proof} \label{proof-sec}

McKenzie's well-known 4-basis for lattices consists of the following
equations.
\begin{align}
x \J (y \M (x \M z)) & = x \tag{L1} \\
x \M (y \J (x \J z)) & = x \tag{L2} \\
((y \M x) \J (x \M z)) \J x & = x \tag{L3} \\
((y \J x) \M (x \J z)) \M x & = x \tag{L4}
\end{align}
The following (machine-oriented) proof is a derivation of
\{L1,L2,L3,L4\} from (A1).
If $i$ stands for the equation $u_i=v_i$, and $j$ stands for $u_j=v_j$,
then [$i\rightarrow j$] justifies an equation $u=v$ obtained in the
following way.  Take a subterm $s$ of $u_j$ such that $u_i$ and $s$ 
are unifiable by substitutions $f$ and $g$; that is, $f(u_i)$ and $g(s)$
are identical.
Let $u$ be the term $g(u_j)$ in which an occurrence of $g(s)$
was replaced by $f(v_i)$, and let $v$ be the term $g(v_j)$.
{\small
\begin{tabbing} 0000\=00000\=\kill
 1 \> $(((y\J x)\M x)\J (((z\M  (x\J x))\J (u\M x))\M v))\M  (w\J ((s\J x)\M  (x\J t)))=x$ \` [A1]\\
 2 \> $(((x\J y)\M y)\J (y\J y))\M  (z\J ((u\J y)\M  (y\J v)))=y$ \` [1 $\rightarrow$ 1] \\
 3 \> $(((x\J (y\J y))\M  (y\J y))\J ((y\J y)\J (y\J y)))\M  (z\J y)=y\J y$ \` [2 $\rightarrow$ 2] \\
 4 \> $(((x\J y)\M y)\J (((y\J y)\J (z\M y))\M u))\M  (v\J ((w\J y)\M  (y\J t)))=y$ \` [3 $\rightarrow$ 1] \\
 5 \> $(((x\J (((y\J y)\J (z\M y))\M u))\M  (((y\J y)\J (z\M y))\M u))\J ((((y\J y)\J (z\M y))\M u)\J$ \\ \>\> $(((y\J y)\J (z\M y))\M u)))\M  (v\J y)= ((y\J y)\J (z\M y))\M u$ \` [4 $\rightarrow$ 2] \\
 6 \> $(((x\J y)\M y)\J (((((y\J y)\J (z\M y))\M u)\J (v\M y))\M w))\M  (t\J ((s\J y)\M  (y\J r)))=y$ \` [5 $\rightarrow$ 1] \\
 7 \> $(((x\J y)\M y)\J (z\M y))\M  (u\J ((v\J y)\M  (y\J w)))=y$ \` [6 $\rightarrow$ 6] \\
 8 \> $(((x\J (y\M z))\M  (y\M z))\J (u\M  (y\M z)))\M  (v\J z)=y\M z$ \` [7 $\rightarrow$ 7] \\
 9 \> $(((x\J y)\M y)\J (((z\M y)\J (u\M y))\M v))\M  (w\J ((t\J y)\M  (y\J s)))=y$ \` [8 $\rightarrow$ 1] \\
 10 \> $(((x\J y)\M y)\J y)\M  (z\J ((u\J y)\M  (y\J v)))=y$ \` [9 $\rightarrow$ 9] \\
 11 \> $(((x\J y)\M y)\J (((z\M y)\J (u\M y))\M v))\M  (w\J y)=y$ \` [10 $\rightarrow$ 9] \\
 12 \> $(((x\J y)\M y)\J (z\M y))\M  (u\J y)=y$ \` [10 $\rightarrow$ 7] \\
 13 \> $(((x\J y)\M y)\J (y\J y))\M  (z\J y)=y$ \` [10 $\rightarrow$ 2] \\
 14 \> $(x\J (y\M  (x\J x)))\M  (z\J ((u\J (x\J x))\M  ((x\J x)\J v)))=x\J x$ \` [13 $\rightarrow$ 7] \\
 15 \> $(x\J x)\M  (y\J ((z\J (x\J x))\M  ((x\J x)\J u)))=x\J x$ \` [13 $\rightarrow$ 14] \\
 16 \> $(((x\J y)\M y)\J ((z\M y)\J (z\M y)))\M  (u\J y)=y$ \` [15 $\rightarrow$ 11] \\
 17 \> $((x\M y)\J (x\M y))\M  (z\J y)= (x\M y)\J (x\M y)$ \` [16 $\rightarrow$ 15] \\
 18 \> $((x\J y)\M y)\J ((x\J y)\M y)=y$ \` [12 $\rightarrow$ 17] \\
 19 \> $(x\M y)\M  (z\J y)=x\M y$ \` [18 $\rightarrow$ 8] \\
 20 \> $x\M  (y\J x)=x$ \` [18 $\rightarrow$ 12] \\
 21 \> $x\M  (y\J ((z\J x)\M  (x\J u)))=x$ \` [18 $\rightarrow$ 7] \\
 22 \> $(((x\J y)\M y)\J ((z\M y)\J (u\M y)))\M  (v\J y)=y$ \` [20 $\rightarrow$ 11] \\
 23 \> $((x\M y)\J (z\M y))\M  (u\J y)= (x\M y)\J (z\M y)$ \` [22 $\rightarrow$ 21] \\
 24 \> $((x\J y)\M y)\J (z\M y)=y$ \` [12 $\rightarrow$ 23] \\
 25 \> $x\M  (x\J y)=x$ \` [24 $\rightarrow$ 21] \\
 26 \> $((x\J y)\M y)\J ((z\M y)\J (u\M y))=y$ \` [22 $\rightarrow$ 25] \\
 27 \> $((x\J y)\M y)\J y=y$ \` [10 $\rightarrow$ 25] \\
 28 \> $x\M  (y\J (x\J z))=x\M  (x\J z)$ \` [25 $\rightarrow$ 19] \\
 29 \> $(x\M x)\J x=x$ \` [27 $\rightarrow$ 27] \\
 30 \> $x\M  (y\J (x\M  (x\J z)))=x$ \` [27 $\rightarrow$ 21] \\
 31 \> $x\M  ((y\J x)\M  (x\J z))=x$ \` [27 $\rightarrow$ 21] \\
 32 \> $x\M x=x$ \` [27 $\rightarrow$ 20] \\
 33 \> $x\M  (y\J (x\J z))=x$ \` (L2) \ \ \ [25 $\rightarrow$ 28] \\
 34 \> $x\J x=x$ \` [32 $\rightarrow$ 29] \\
 35 \> $(x\J y)\M y=y$ \` [24 $\rightarrow$ 34] \\
 36 \> $x\M  (y\M x)=y\M x$ \` [24 $\rightarrow$ 35] \\
 37 \> $(x\J (((y\M x)\J (z\M x))\M u))\M  (v\J x)=x$ \` [35 $\rightarrow$ 11] \\
 38 \> $((x\J y)\M  (y\J z))\M y=y\M  ((x\J y)\M  (y\J z))$ \` [31 $\rightarrow$ 36] \\
 39 \> $(x\J (x\M y))\M  (z\J x)=x$ \` [24 $\rightarrow$ 37] \\
 40 \> $x\J (((y\M x)\J (z\M x))\M u)=x$ \` [25 $\rightarrow$ 37] \\
 41 \> $((x\J y)\M  (y\J z))\M y=y$ \` (L4) \ \ \ [31 $\rightarrow$ 38] \\
 42 \> $(((x\M y)\J (z\M y))\J (((x\M y)\J (z\M y))\M u))\M y= (x\M y)\J (z\M y)$ \` [26 $\rightarrow$ 39] \\
 43 \> $(x\M y)\M x=x\M y$ \` [39 $\rightarrow$ 31] \\
 44 \> $x\J (((y\M x)\J x)\M z)=x$ \` [35 $\rightarrow$ 40] \\
 45 \> $x\J ((y\M x)\M z)=x$ \` [40 $\rightarrow$ 40] \\
 46 \> $((x\M y)\J y)\M y= (x\M y)\J y$ \` [44 $\rightarrow$ 30] \\
 47 \> $x\J (y\M  (z\M x))=x$ \` [36 $\rightarrow$ 45] \\
 48 \> $(x\M y)\J y=y$ \` [35 $\rightarrow$ 46] \\
 49 \> $x\J (y\M  (x\M z))=x$ \` (L1) \ \ \ [43 $\rightarrow$ 47] \\
 50 \> $((x\M y)\J (z\M y))\J y=y$ \` [42 $\rightarrow$ 48] \\
 51 \> $((x\M y)\J (y\M z))\J y=y$ \` (L3) \ \ \ [43 $\rightarrow$ 50] \\
 
\end{tabbing}
} %
\noindent
Lines \{49,33,51,41\} are equations \{L1,L2,L3,L4\}, respectively.
A similar \otter\ proof shows that (A2) is a single identity.

The Web page
\begin{quote} 
\texttt{www.mcs.anl.gov/\~{}mccune/papers/ltsax}
\end{quote}
contains several files associated with this note, including
previously known single identities,
\otter\ input files that produced the proofs,
and other supporting material.

\bibliographystyle{plain}

\end{document}